# Conditional probability of actually detecting a financial fraud – a neutrosophic extension to Benford's law

Sukanto Bhattacharya
Alaska Pacific University, USA

Kuldeep Kumar
Bond University, Australia

Florentin Smarandache
University of New Mexico, USA

**Abstract**
___________________________________________________________
This study actually draws from and builds on an earlier paper (Kumar and Bhattacharya, 2002). Here we have basically added a *neutrosophic dimension* to the problem of determining the *conditional probability* that a financial fraud has been actually committed, given that no Type I error occurred while rejecting the null hypothesis **H₀:** The observed first-digit frequencies approximate a Benford distribution; and accepting the alternative hypothesis **H₁:** The observed first-digit frequencies do not approximate a Benford distribution. We have also suggested a conceptual model to implement such a neutrosophic fraud detection system.
___________________________________________________________

**Re-visiting the problem of testing for manipulation in accounting data**

In an earlier paper (Kumar and Bhattacharya, 2002), we had proposed a Monte Carlo adaptation of Benford's law. There has been some research already on the application of Benford's law to financial fraud detection (Carslaw, 1988 and Busta and Weinberg 1998). However, most of the practical work in this regard has been concentrated in detecting the first digit frequencies from the account balances selected on basis of some known audit sampling method and then directly comparing the result with the expected Benford frequencies (Raimi, 1976 and Hill, 1998). We have voiced our reservations about this technique in so far as that the Benford frequencies are necessarily *steady state frequencies* and may not therefore be truly reflected in the sample frequencies. As samples are always of finite sizes, it is therefore perhaps not entirely fair to arrive at any conclusion on the basis of such a direct comparison, as the sample frequencies won't be steady state frequencies.

However, if we draw digits randomly using the *inverse transformation technique* from within random number ranges derived from a cumulative probability distribution function based on the Benford frequencies then the problem boils down to running a *goodness of fit* kind of test to identify any significant difference between observed and simulated first-digit frequencies. This test may be conducted using a known sampling distribution like for example the *Pearson's $\chi^2$ distribution*. The random number ranges for the Monte Carlo simulation are to be drawn from a cumulative probability distribution function based on the following Benford probabilities given in Table I.

**Table I**

| First Significant Digit | 1 | 2 | 3 | 4 | 5 | 6 | 7 | 8 | 9 |
|---|---|---|---|---|---|---|---|---|---|
| Benford Probability | 0.301 | 0.176 | 0.125 | 0.097 | 0.079 | 0.067 | 0.058 | 0.051 | 0.046 |

The first-digit probabilities can be best approximated mathematically by the log-based formula Benford derived: **P (First significant digit = d) = $\log_{10}$ [1 + (1/d)]** (Benford, 1938).

**Computational Algorithm:**

1. Define a finite sample size n and draw a sample from the relevant account balances using a suitable audit sampling procedure

2. Perform a continuous Monte Carlo run of length $\lambda^* \approx (1/2\varepsilon)^{2/3}$ grouped in epochs of size n using a customized MS-Excel spreadsheet. Derivation of $\lambda^*$ and other statistical issues have been discussed in detail in our earlier paper (Kumar and Bhattacharya, 2002)

3. Test for significant difference in sample frequencies between the first digits observed in the sample and those generated by the Monte Carlo simulation by using a "goodness of fit" test using the $\chi^2$ distribution. The null and alternative hypotheses are as follows:

**H₀:** The observed first digit frequencies approximate a Benford distribution

**H₁:** The observed first digit frequencies do not approximate a Benford distribution

This statistical test will not reveal whether or not a fraud has actually been committed. All it does is establishing at a desired level of confidence, that the accounting data has been manipulated (if $H_0$ is rejected).

However, given that $H_1$ is accepted and $H_0$ is rejected, it could imply any of the following events:

I. There is no manipulation - Type I error has occurred i.e. $H_0$ rejected when true.

II. There is manipulation *and* such manipulation *is definitely* fraudulent.

III. There is manipulation *and* such manipulation *may or may not be* fraudulent.

IV. There is manipulation *and* such manipulation *is definitely not* fraudulent.

**Neutrosophic extension**

Neutrosophic probabilities are a generalization of classical and fuzzy probabilities and cover those events that involve some degree of indeterminacy. It provides a better approach to quantifying uncertainty than classical or even fuzzy probability theory. Neutrosophic probability theory uses a subset-approximation for truth-value as well as indeterminacy and falsity values. Also, this approach makes a distinction between "relative true event" and "absolute true event" the former being true in only some probability sub-spaces while the latter being true in all probability sub-spaces. Similarly, events that are false in only some probability sub-spaces are classified as "relative false events" while events that are false in all probability sub-spaces are classified as "absolute false events". Again, the events that may be hard to classify as either 'true' or 'false' in some probability sub-spaces are classified as "relative indeterminate events" while events that bear this characteristic over all probability sub-spaces are classified as "absolute indeterminate events". (Smarandache, 2001)

While in classical probability **n_sup $\leq$ 1**, in neutrosophic probability **n_sup $\leq$ 3$^+$** where n_sup is the upper bound of the probability space. In cases where the truth and falsity components are complimentary, i.e. there is no indeterminacy, the components sum to unity and neutrosophic probability is reduced to classical probability as in the tossing of a fair coin or the drawing of a card from a well-shuffled deck.

Coming back to our original problem of financial fraud detection, let E be the event whereby a Type I error has occurred and F be the event whereby a fraud is actually detected. Then the *conditional neutrosophic probability* **NP (F | E$^c$)** is defined over a probability space consisting of a triple of sets (T, I, U). Here, T, I and U are probability sub-spaces wherein event F is t% true, i% indeterminate and u% untrue respectively, given that no Type I error occurred.

The sub-space T within which t varies may be determined by factors such as past records of fraud in the organization, propensity to commit fraud by the employees concerned, and effectiveness of internal control systems. On the other hand, the sub-space U within which u varies may be determined by factors like personal track records of the employees in question, the position enjoyed and the remuneration drawn by those employees. For example, if the magnitude of the embezzled amount is deemed too frivolous with respect to the position and remuneration of the employees involved. The sub-space I within which i varies is most likely to be determined by the mutual inconsistency in the circumstantial evidence (Zadeh, 1976) that might arise out of the effects of some of the factors determining T and U.  For example, if an employee is for some reason really irked with the organization, then he or she may be inclined to commit fraud not so much to further his or her own interests as to harm the interests of the organization, although the act of actually committing the suspected fraud may in this case overtly appear inconsistent with the organizational status and remuneration enjoyed by that person.

**A conceptual model to implement the neutrosophic fraud detection system**

Modern technology has armed the investigative accountants with tools and techniques not only to track down the perpetrators of fraud more efficiently than was possible in the absence of those technologies but also to carry out a multi-faceted analytical inquiry into the nature of financial frauds and their perpetrators.

We propose classifying financial frauds in the modern corporate scenario using a systematic, multi-level categorization. The simplest one would of course be a two-level classification where one classifier is in terms of the *nature of fraudulent manipulation* and the other is in terms of *involvement of the perpetrators* e.g. fraud by an individual, a group of isolated individuals or a collusive fraud. A sub-classification within this second categorization could be based on whether the group (isolated or collusive) consists of hierarchical positions or horizontal positions in the organizational structure. Then what may come up with is the following two-dimensional *manipulation-involvement* or *MI*-matrix:

Chart 1:
Increasing complexity of involvement (i)

|  | Manipulation | | | | |
|---|---|---|---|---|---|
| Involvement | $\alpha_{11}$ | $\alpha_{12}$ | $\alpha_{13}$ | $\alpha_{14}$ | $\alpha_{15}$ |
|  | $\alpha_{21}$ | $\alpha_{22}$ | $\alpha_{23}$ | $\alpha_{24}$ | $\alpha_{25}$ |
|  | $\alpha_{31}$ | $\alpha_{32}$ | $\alpha_{33}$ | $\alpha_{34}$ | $\alpha_{35}$ |
|  | $\alpha_{41}$ | $\alpha_{42}$ | $\alpha_{43}$ | $\alpha_{44}$ | $\alpha_{45}$ |
|  | $\alpha_{51}$ | $\alpha_{52}$ | $\alpha_{53}$ | $\alpha_{54}$ | $\alpha_{55}$ |
|  | $\alpha_{61}$ | $\alpha_{62}$ | $\alpha_{63}$ | $\alpha_{64}$ | $\alpha_{65}$ |
|  | $\alpha_{71}$ | $\alpha_{72}$ | $\alpha_{73}$ | $\alpha_{74}$ | $\alpha_{75}$ |
|  | $\alpha_{81}$ | $\alpha_{82}$ | $\alpha_{83}$ | $\alpha_{84}$ | $\alpha_{85}$ |

Increasing complexity of manipulation (j)

The row and column elements in this simple *MI*-matrix may be denoted as follows:

$\alpha_{ij}$:

i = 8: Single fraud perpetrated by single individual

i = 7: Multiple frauds perpetrated by single individual

i = 6: Single fraud perpetrated by group of isolated individuals

i = 5: Multiple frauds perpetrated by group of isolated individuals

i = 4: Single fraud perpetrated by a horizontally collusive group

i = 3: Multiple frauds perpetrated by a horizontally collusive group

i = 2: Single fraud perpetrated by a vertically collusive group

i = 1: Multiple frauds perpetrated by a vertically collusive group

j = 5: Combination of 1, 2, 3 and 4

j = 4: Suppression or destruction of key transaction records

j = 3: Misrepresentation of the fundamental nature of transaction

j = 2: Falsification of transaction date/amount/particulars (double-entry basis)

j = 1: Falsification of transaction date/amount/particulars (single-entry basis)

Once a particular case has been objectively classified on the *MI*-matrix, the investigative accountant may start looking for incriminating evidence in the financial records *from the right perspective*. For example, the search perspective will definitely differ if there is fundamental alteration in the nature of a transaction as compared to a simple erring journal entry. The perspective will also differ if only a single individual is involved rather than a collusive group.

The two-dimensional *MI*-matrix is the simplest form of systematic multi-level classification which certainly may be conceptually expanded to include more than two levels – e.g. a three-dimensional *MI*-matrix could possibly incorporate the aspect of *fraud potentiality* in terms of factors like personal track records of the employees in question, the position enjoyed by them in the organization, their remuneration and entitlements at the time of the fraud etc. It would therefore be an ideal computational set-up for implementation of a neutrosophic system. The conceptual fraud classification scheme we proposed here may thus be effectively combined with the rules of neutrosophic probability into developing a handy *forensic accounting expert system* for the future!

**Conclusion**

No doubt then that the *theory of neutrosophic probability* opens up a new vista of analytical reasoning for the techno-savvy forensic accountant. In this paper, we have only posit that a combination of statistical testing of audit samples based on Benford's law together with a neutrosophic reasoning could help the forensic accountant in getting a better fix on the quantitative possibility of actually detecting a financial fraud. This is an emerging science and thus holds a vast potential of future research endeavours the ultimate objective of which will be to actually come up with a reliable, comprehensive computational methodology to actually track down financial frauds with a very low failure rate. We believe our present work is just an initial step towards that ultimate destination.